\documentclass[11pt,a4paper]{article}
\usepackage{authblk}

\usepackage{graphicx}

\usepackage{amssymb}
\usepackage{amsmath}
\usepackage{amsfonts}
\usepackage{color}

\newcommand{\T}{{\mathbb T}}
\newcommand{\R}{{\mathbb R}}
\newcommand{\C}{{\mathbb C}}

\newcommand{\Proj}{{\mathbb P}}

\newtheorem{theorem}{Theorem}[section]
\newtheorem{proposition}[theorem]{Proposition}

\newenvironment{Proof}{\removelastskip \noindent {\bf Proof~:} } { \hspace*{\fill} {\bf q.e.d.} \medskip \noindent}

\title{Poisson and symplectic reductions of 4-DOF isotropic
oscillators. The van der Waals system as benchmark}
\author[1]{F. Crespo}
\author[1]{G. D\'{\i}az--Toca}
\author[1]{S. Ferrer}
\author[2]{M. Lara}
\affil[1]{{\small Dpto de Matem\'atica Aplicada, Universidad de Murcia, 30071 Espinardo,  Spain}}
\affil[2]{{\small Observatorio de la Armada, 11110 San Fernando, Spain}}

\begin{document}
\maketitle

\begin{abstract}
This paper is devoted to studying Hamiltonian oscillators in 1:1:1:1 resonance with
symmetries, which include several models of perturbed Keplerian
systems. Normal forms are computed in Poisson and
symplectic formalisms,  by mean of invariants and Lie-transforms respectively. The first procedure relies on the quadratic
invariants associated to the symmetries, and is carried out using
Gr\"obner bases. In the symplectic approach, hinging on
the maximally superintegrable character of the isotropic oscillator,
the normal form is computed  {\it a la} Delaunay, using a
generalization of those variables for 4-DOF systems. Due to the symmetries of
the system, isolated as well as circles of stationary
points and invariant tori should be expected. These solutions manifest  themselves rather differently in 
both approaches, due to the constraints among the invariants versus the singularities associated to the Delaunay chart.\par

Taking the generalized van der Waals family as a benchmark, the explicit expression of the Delaunay normalized Hamiltonian  up to the second order is presented, showing that it may be extended to higher orders in a straightforward way. The search for the relative equilibria is used for comparison of their main features of both treatments.  The pros and cons are given in detail for some values of the parameter and the
integrals.
\end{abstract}


\section{Introduction}\label{intro}

The use of computer algebra systems for normal forms computations is
considered at present a routine operation. As a general reference see
e.g. Sanders {\it et al.} \cite{sanders2007} and Meyer {\it et al.} \cite{Meyer09}.  Nevertheless when we
deal with special classes of differential equations, like Poisson or
Hamiltonians systems which is our case, it is advisable to employ
specific transformations as well as tailored variables for those
problems \cite{Meyer09}, mostly connected with the symmetries that
those systems might possess. More precisely we are interested in
perturbed isotropic oscillators in four dimensions (or perturbed
harmonic oscillators in 1:1:1:1 resonance). For those
Hamiltonian systems,  Sanders {\sl et al.} explain that the 1:...:1
resonance is one of the more complicated, due to large number
of terms in the normal form. This proves to be a key
issue in computations of higher order approximations, which is needed in bifurcation analysis for some values of the parameters. For this reason,
any strategy developed to reduce the algebra involved in the normal 
form process, as well as in the subsequent analysis built on it, like
relative equilibria and their bifurcations, is something really
desirable. As we will see, we have to confront with systems of polynomial equations with parameters, which represent a real challenge, even with computer algebra assistance.\par

Continuing previous work \cite{cushman1999,Ferrer2002,Ferrer2000} \cite{Ferrer2002} in the 1:1:1 resonance,  and also \cite{egea2007a,egea2007b,diaz2009} in relation to 1:1:1:1 resonance,  we consider in $\mathbb{R}^{4}\times\mathbb{R}^{4} $,
the symplectic form $w=d\mathbf{Q\wedge
dq}$ and a Hamiltonian function
\begin{equation}\label{hamiltoniano0}
\mathcal{H}=\mathcal{H}_2+ \varepsilon \mathcal{P}(q,Q;\beta)
\end{equation}
where
\begin{equation}\label{h2}
\mathcal{H}_{2}=\frac{1}{2}(Q_{1}^{2}+Q_{2}^{2}+Q_{3}^{2}+Q_{4}^{2})+\frac
{1}{2}\,\omega^2\,(q_{1}^{2}+q_{2}^{2}+q_{3}^{2}+q_{4}^{2})
\end{equation}
defines the isotropic oscillator, with $\omega$ a positive constant
and $\varepsilon$ is an small parameter. In the first part of the
paper we simplify expressions assuming $\omega=1$. The function
$\mathcal{P}$ is called the perturbation, where $\beta$ is a parameter
vector, which may include also $\varepsilon$. Moreover we consider
that systems defined by Hamiltonian function Eq. (\ref{hamiltoniano0}) have
two first integrals in involution given
by
\begin{align}\label{integrales}
\Xi&= q_{1}Q_{2}-Q_{1}q_{2} + q_{3}Q_{4}-Q_{3}q_{4},\notag \\
L_{1}&=q_{3}Q_{4}-Q_{3}q_{4} - q_{1}Q_{2}+Q_{1}q_{2},
\end{align}
associated to which we have rotational symmetries. We use the same
notation as in \cite{egea2007b}.

Although some of the methods and techniques considered in the paper may
be applied to a large family of Hamiltonian systems defined by
(\ref{hamiltoniano0}),  to give details of those processes we
focus on the uniparametric family of Hamiltonian systems defined by
\begin{equation}\label{Ham1}
\mathcal{H}_{\beta}(Q,q)=\mathcal{H}_{2}+\varepsilon\,\mathcal{H}_{6},
\end{equation}
where
\begin{eqnarray}
&&\hspace{-1cm}\mathcal{H}_{6}(Q,q;\beta)=\left(q_{1}^{2}+q_{2}^{2}+q_{3}
^{2}+q_{4}^{2}\right) \nonumber\\
&&\hspace{1cm}\times \left(  {\beta}^{2}\,{\left(  q_{1}^{2}+q_{2}
^{2}-q_{3}^{2}-q_{4}^{2}\right)  }^{2}+4\, \left(  q_{1}^{2}+q_{2}
^{2}\right)\left(  q_{3}^{2}+q_{4}^{2}\right)  \right)
\end{eqnarray}
where $\beta$ is now a real parameter. For $\beta=1$ we have a central potential, {\it i.e.} an integrable system. When the system restricts to the manifold $\Xi=0$ then it is equivalent
to the model for the hydrogen atom subject
to a generalized van der Waals potential. For $\beta=0$ this system
reduces to the model for the quadratic Zeeman effect. When
$\beta=\sqrt{2}$ we have the Van der Waals system. For this reason we
have named the system as the \textsl{generalized Van der Waals 4-D
oscillator}. A search  (see \cite{diaz2009}) for some  special solutions of the Hamiltonian system defined by the function (\ref{Ham1}) reveals that there are invariant 2-tori associated to the rotational integrals (\ref{integrales}), which include straight-line orbits over the configuration space.

Normalization and reduction using invariants and Lie-transforms are approached in both Poisson and symplectic formalisms. The first procedure (see \cite{egea2007a}, \cite{cushman1997}) relies on the quadratic invariants associated to the symmetries,
and is based on previous work of the authors. Enlarging previous studies on the 1:1:1 resonance \cite{cushman1999,Ferrer2002},  the Hamiltonian (\ref{Ham1}) is put into normal form with respect
to $\mathcal{H}_{2}$. Considering the truncated normal form  a system is obtained that is invariant under the
the $\Bbb{S}^{1}$-actions corresponding to $\mathcal{H}_{2}$, $\Xi$, and $L_1$. These three actions together generate a $\T^3$-action. Reduction with respect to it leads to a one-degree-of-freedom system. \par

Hinging on the maximally superintegrable character of the isotropic oscillator, the symplectic reduction is carried out {\sl a la} Delaunay using a generalization of those variables to 4-DOF recently proposed in \cite{Ferrer2010}. Based on the normalized equations, our studies on this
system focus on relative equilibria and their
evolution with the parameter $\beta$. Although first order
normalization (averaging) is quite generic, the fact that
the family includes some integrable cases, makes it necessary to perform higher
order normalization for the stability analysis of those cases.

The paper is organized as follows. Section \ref{sec:normalization} summarizes part of the work done in \cite{diaz2010}, the three steps of the Poisson reduction are presented with the invariant functions involved, as well as the corresponding reduced Hamiltonian of our model at each step. In Section \ref{sec:simplecticas} we proceed likewise carrying out the reduction, although in different order, by three symplectic transformations: Projective Euler and  Andoyer as well as 4-D Delaunay. They give again, but now in the symplectic frame, the 1-DOF thrice reduced Hamiltonian system in the open domain where those charts are defined. Once carried out the toral reduction by the symplectic variables in Section \ref{sec:simplecticas}, Section \ref{sec:first} focuses on the computation of the normal form by Lie transform which, taking advantage of the periodic character of the unperturbed flow, which allows it to solve the homological equation just by quadratures. Finally Section \ref{sec:periodic} gathers the analysis on relative equilibria which enlarges the results presented in previous papers, in particular the search for periodic orbits.

\section{Poisson approach. Normalization and reduction}
\label{sec:normalization}
There is a large literature which develops both theoretical and
computational aspects of normal forms, see for example \cite{Broer03,Mayer2004,sanders2007,Meyer09}. In this section we briefly describe
several normal forms for first order of system (\ref{Ham1}) for
$\omega=1$ thoroughly presented in \cite{diaz2010}. In such a work,
we used Gr\"{o}bner Bases as main computational tool and, as
software, the computer algebra system \texttt{Maple}.

A constructive geometric reduction in
stages, including regular and singular reductions, is performed in \cite{diaz2010} and only the final steps are presented here.
In the Poisson approach the first reduction is done with respect to the $H_2$
symmetry. This is a regular reduction and the reduced phase space is
homeomorphic to $\mathbb{CP}^{3}$ (see \cite{cushman1997}). Then we
carry out a second reduction whose resulting orbit space is
stratified: four dimensional isomorphic to $\Bbb{S}^2 \times
\Bbb{S}^2$, and two singular strata isomorphic to $\Bbb{S}^2$.
Finally we make a third reduction, with the integral $L_{1}$,
reducing the system to a 1-DOF system on the \textsl{thrice reduced
phase space}. Depending on the relative value of the integrals, they
are isomorphic or homomorphic to a 2-sphere, containing one or two
singular points. Besides there are singular reduced phase
spaces consisting of a single point (see Figure 2 and
\cite{egea2007b}).

\subsection{Normalization with respect to the oscillator symmetry
$X_{H_2}$}

Given the action associated to the uniparametric group defined by
$X_{H_2}$, by using canonical complex variables, (see
\cite{egea2007b} for details) it follows that the algebra of
polynomial invariants under that action is generated by
\begin{equation}
\begin{array}{lllll}
&\pi_i=Q_i^2 + q_i^2, &\;i=1,2,3,4& & \\
&\pi_5=Q_1\,Q_2 + q_1\,q_2, &\;\pi_6 =Q_1\,Q_3 + q_1\,q_3, &\;\pi_7
=Q_1\,Q_4 + q_1\,q_4,&\\
&\pi_8 =Q_2\,Q_3 + q_2\,q_3,&
\;\pi_9 = Q_2\,Q_4 + q_2\,q_4,&\;\pi_{10}=Q_3\,Q_4 + q_3\,q_4, &\\
&\pi_{11}= q_1\,Q_2-Q_1\,q_2,&\;\pi_{12}= q_1\,Q_3-Q_1\,q_3 ,&
\;\pi_{13}= q_1\,Q_4-Q_1\,q_4, &\\
&\pi_{14}= q_2\,Q_3 -Q_2\,q_3, &\;\pi_{15}= q_2\,Q_4-Q_2\,q_4,&\;
\pi_{16}= q_3\,Q_4-Q_3\,q_4.
\end{array}
\end{equation}
The $X_{H_2}$ normal form up to first order in $\varepsilon$ is expressed in those invariants as
\begin{equation}
\overline{\mathcal{H}}=\mathcal{H}_{2}+\varepsilon\overline{\mathcal{H}}_{6}%
\end{equation}
where
$\mathcal{H}_{2}=\left(\pi_{1}+\pi_{2}+\pi_{3}+\pi_{4}\right)/2=n$ and
\begin{eqnarray}
&&\hspace{-1.2cm}\overline{\mathcal{H}}_{6} = 
\frac{1}{2}\left[
(1-4\beta^{2})n\,(\pi_{15}^{2} +\pi_{14}^{2}+\pi_{13}^{2}+\pi_{12}^{2})\right.\notag\\
&&+2(\beta^{2}-1)(\pi_{11}^{2}(\pi_{4}+\pi_{3})-\pi_{16}^{2}(\pi_{3}+\pi_{4})) +\beta^{2} n\,(5n^{2}-3\pi_{11}^{2})\\
&&\left.+5(1-\beta^{2})n\,(\pi_{9}^{2}+\pi_{8}^{2}+\pi_{7}^{2}+\pi_{6}^{2})+(\beta^{2}-4)n\,\pi_{16}^{2}\right]\notag
\end{eqnarray}
The reduction is now performed using the orbit map
$$\rho_{\pi}:\R^{8} \rightarrow \R^{16};\,\, (q,Q)\rightarrow (\pi_1,
\cdots ,\pi_{16})\, .$$
The image of this map is the orbit space for the $H_2$-action, the
images of the level surfaces $H_2(q,Q)=n$ under $\rho_{\pi}$ are the
reduced phase spaces which are isomorphic to $\C\Proj^{3}$. The
normalized Hamiltonian can be expressed in
invariants and therefore naturally lifts to a function on $\R^{16}$,
which, on the reduced phase spaces, restricts to the reduced
Hamiltonian.

However, in the following we will not use the invariants $\pi_i$;
instead, we rely on $(K_i,L_j,J_k)$ invariants introduced in Egea
\cite{egea2007a} by the following change
of coordinates,
\[
\begin{array}{ll}
 H_2 =\frac{1}{2}(\pi_1 +\pi_2 +\pi_3 +\pi_4),  & \qquad J_1
=\frac{1}{2}(\pi_1 -\pi_2 -\pi_3 +\pi_4),\\[1ex]
 \hspace{0.2cm}\Xi = \pi_{16} +\pi_{11},  & \qquad   J_2
=\frac{1}{2}(\pi_1 -\pi_2 +\pi_3 -\pi_4),\\[1ex]
 K_1 =\frac{1}{2}(-\pi_1 -\pi_2 +\pi_3 +\pi_4), &\qquad  J_3 = \pi_8
+\pi_7, \\[1ex]
 K_2 = \pi_8 -\pi_7, & \qquad  J_4 = \pi_5 +\pi_{10}, \\[1ex]
 K_3 = -\pi_6 -\pi_9, &\qquad  J_5 = \pi_5-\pi_{10},\\[1ex]
 L_1 = \pi_{16}-\pi_{11}, & \qquad  J_6 = \pi_6 -\pi_9 ,\\[1ex]
 L_2 = \pi_{12} +\pi_{15}, &\qquad   J_7 = \pi_{12} -\pi_{15},\\[1ex]
 L_3 = \pi_{14} -\pi_{13}, &\qquad   J_8 = \pi_{14} +\pi_{13}.
\end{array}
\]
The {\sl  first integrals} (see Eq. \ref{integrales}) are now among
the invariants defining the image.

The first order normal form  in these invariants is
\begin{align}\label{primerreducido}
\overline{\mathcal{H}}_{\Xi}&=\frac{1}{2}\Big[n\,\Big(5\,K{_{2}}^{2}+5\,K{_{3}}
^{2}+2\,L{_{1}}^{2}+L{_{2}}^{2}+L{_{3}}^{2}  +
{\beta}^{2}\,(5\,K{_{1}}^{2}+L{_{2}}^{2} +L{_{3}}^{2})\Big)\nonumber\\
&\hspace{0.8cm}-\left( (4+{\beta}^{2})\,(K_{2}\,L_{2}+K_{3}\,L_{3}) +
(2+3\,{\beta}^{2})K_{1} L_{1}\right)\,\xi\Big]\, .
\end{align}
The reduction of the $H_2$ action may now be performed through the
orbit map
$$
\rho_{K,L,J}:\R^8 \rightarrow \R^{16}; (q,Q)\rightarrow (H_2, \cdots
,J_8) \; .
$$
Note that on the orbit space we have the reduced symmetries due to
the reduced actions given by the reduced flows of $X_{\Xi}$ and
$X_{L_1}$.

\subsection{Toroidal reduction over $\C\Proj^{3}$ with respect to the
rotational symmetry $\Xi$}
\begin{figure}
\begin{center}
\includegraphics[width=250pt]{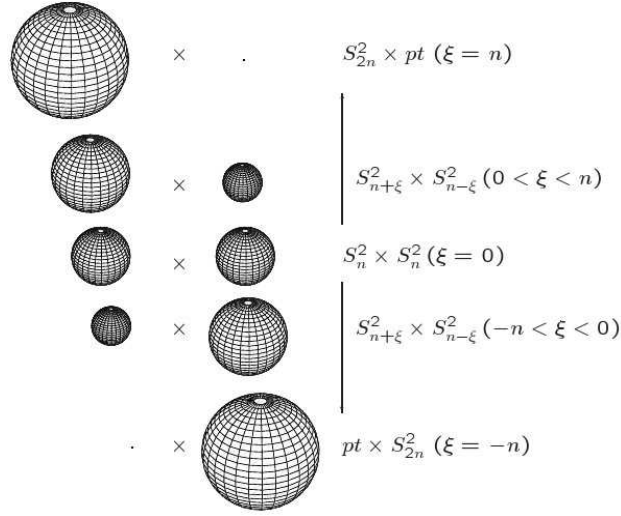}
\end{center}
\caption{Double reduced spaces $S_{n+\xi}^2\times S_{n-\xi}^2$ for
different values of the integral $\xi$}
\label{fig:DoubleReducedSpaces}
\end{figure}

Let $\rho$ be the $S^1$-action generated by the Poisson flow of
$\Xi$ over $\C\Proj^3$. The functions
\[H_2,\quad \Xi, \quad L_1,\quad L_2,\quad L_3,\quad K_1,\quad
K_2,\quad K_3,\]
are $\rho$-invariants over $\C\Proj^3$. This, in turn, leads us to
the orbit mapping
\[\rho_2 : \R^{16} \rightarrow \R^{8}; \,\,(\pi_1 ,\cdots
,\pi_{16})\rightarrow (K_1,K_2,K_3,L_1,L_2,L_3,H_2,\Xi).\]
The orbit space $\rho_2(\C\Proj^3)$ is defined as a six dimensional
algebraic variety in $\R^{8}$ by the following relations
\begin{equation}
K_1^2 +K_2^2 +K_3^2+L_1^2 +L_2^2 +L_3^2= H_2^2 + \Xi^2\; ,\quad
K_1 L_1 +K_2 L_2 +K_3 L_3=H_2\, \Xi\;
\end{equation}
The reduced phase spaces are obtained by setting $ H_2=n$,  $\Xi=\xi
$ and then the second reduced space is isomorphic to $S^2_{n+\xi}\times
S^2_{n-\xi}$ (see Fig.~\ref{fig:DoubleReducedSpaces}).
When  $\xi = 0$, it corresponds to the first reduced space of
Keplerian
systems by the energy.
The second reduced Hamiltonian up to first order, modulo a constant
takes the form
\begin{eqnarray*}
&\overline{\mathcal{H}}_{\Xi}=\frac{1}{2}
\left[n\,\left(5\,K{_{2}}^{2}+5\,K{_{3}}^{2}+2\,L{_{1}}^{2}+L{_{2}}^{2}+L{_{3}}^{2}+{\beta}^{2}\,(5\,K{_{1}}^{2}+L{_{2}}^{2}
+L{_{3}}^{2})\,\right)\right.\\
&\left.-\left( (4+{\beta}^{2})\,(K_{2}\,L_{2}+K_{3}\,L_{3}) +
(2+3\,{\beta}^{2})K_{1} L_{1}\right)\,\xi\right]
\end{eqnarray*}

\subsection{Reduction by $L_1=l$. Thrice reduced space
$V_{n\,\zeta\,l}$}
To further reduce from $S_{n+\xi}^2\times S_{n-\xi}^2$ to
$V_{n\,\zeta\,l}$, one divides out the $S^{1}$-action, $\rho_2$,
generated by the Poisson flow defined by $L_1=\pi_{16}-\pi_{11}$ over
$S_{n+\xi}^2\times S_{n-\xi}^2$. The 8 invariants for
the $L_1$ action on $\R^{8}$ are
\begin{eqnarray}\label{invariantes3}
&& M = \dfrac{1}{2} (K_2^2 +K_3^2+L_2^2 +L_3^2)\; , \qquad N =
\dfrac{1}{2} (K_2^2 +K_3^2-L_2^2-L_3^2) \; , \nonumber\\
&& Z= K_2 L_2 +K_3 L_3 \; ,  \hspace{2.7cm}S  = K_2 L_3 -K_3 L_2 \;
,\\
&& K = K_1,	\qquad \mathcal{H}_{2},\qquad \Xi, \qquad L_1 \;
.\nonumber
\end{eqnarray}
There are 2 + 3 relations defining the third reduced phase space
\begin{alignat}{1}\label{eo3}
   & K^2 +L_1^2 +2M= n^2 +\xi^2\; ,\quad  M^2-N^2=Z^2+S^2 \;,\notag \\
   &\hspace{1cm} K L_1 + Z=n \xi\; ,\quad
\quad\mathcal{H}_{2}=n,\quad \Xi=\xi, \quad L_1=l \;.
   \end{alignat}
The Poisson structure  of the variables $S$, $K$ and $N$ is given in Table~\ref{poissonmn}.
\begin{table}[htb]
\begin{center}
\begin{tabular}{|c|cccccc|}
\hline
$\{,\}$ & $M$ & $N$ & $Z$ & $S$ & $K$ & $L_1$\\
\hline
$M$ & $0$ & $4KS$ & $0$ & $-4KN$ & $0$ & $0$\\
$N$ & $-4KS$ & $0$ & $-4L_1S$ & $ -4(KM - L_1 Z )$ & $4S$ & $0$\\
$Z$ & $0$ & $4 L_1 S$ & $0$ & $-4 L_1 N$ & $0$ & $0$\\
$S$ & $4 KN$ & $4( K M - L_1 Z ) $ & $4 L_1 N$&$0$ & $-4N$ & $0$\\
$K$ & $0$ & $-4S$ & $0$ & $4N$ & $0$ & $0$\\
$L_1$& $0$ & $0$ & $0$ & $0$ & $0$ & $0$\\
\hline
\end{tabular}\caption{Poisson structure in $(M,N,Z,S,K,L_1)$ invariants}
\label{poissonmn}
\end{center}
\end{table}  \\
For more details see \cite{diaz2010}. Then, we have defined the orbit map
\[
\rho_2: \R^{6} \rightarrow \R^{6};
(K_1,K_2,K_3,L_1,L_2,L_3)\rightarrow (M,N,Z,S,K,L_1).
\]
\begin{figure}[htb]
\begin{center}
\includegraphics[width=250pt]{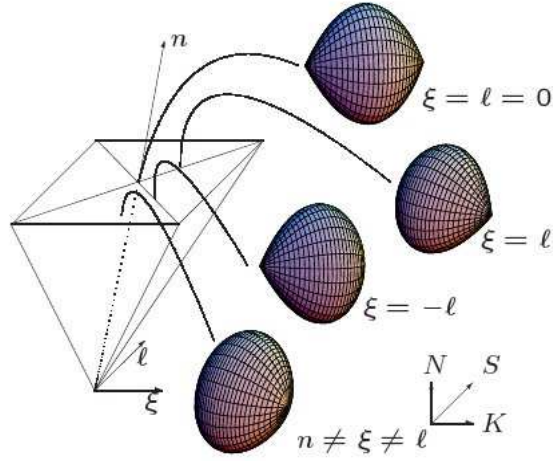}
\end{center}
\caption{\small{Thrice reduced space over the space of integrals. The
vector $(K,N,S)$
represents the coordinates. The axis of symmetry of the reduced space is the $K$ direction.}}
\label{fig:ThriceReducedSpace}
\end{figure}
When we fix a value of  $L_1=l$, relations (\ref{eo3}) define the
{\sl thrice reduced space}
\begin{eqnarray*}\label{eredmn}
&&V_{n\,\xi\,l} = \{\,(K,S,N) \,\,|\,\,  4N^2+4S^2=f(K),\\[1.2ex]
&&\hspace{2cm} f(K)=((n+\xi)^2-(K+l)^2)((n-\xi)^2-(K-l)^2)\,\}
\end{eqnarray*}\noindent
which is a surface of revolution, obtained by rotating $\sqrt{f(K)}$
around the axis $K$, . Thus the shape of the
reduced phase space, given in Fig.~\ref{fig:ThriceReducedSpace}, is determined by the positive part of $f(K)$.
Since
\[
f(K)=(K+n+\xi+l)(K-n-\xi+l)(K-n+\xi-l)(K+n-\xi-l),
\]
the roots are
\[
k_1=-l-n-\xi, \quad k_2=l+n-\xi, \quad k_3=l-n+\xi, \quad
k_4=-l+n+\xi \; .
\]
So $f(K)$ is positive (or zero) in the subsequent intervals of $K$:
\begin{align}
l<\xi\,,\,-l<\xi & \quad  k_{1}<k_{3}<k_{2}<k_{4}\quad K\in\lbrack
k_{3},k_{2}]\label{nxil}\notag\\
l>\xi\,,\,-l<\xi & \quad  k_{1}<k_{3}<k_{4}<k_{2}\quad K\in\lbrack
k_{3},k_{4}]\notag\\
l<\xi\,,\,-l>\xi & \quad k_{3}<k_{1}<k_{2}<k_{4}\quad K\in\lbrack
k_{1},k_{2}]\\
l>\xi\,,\,-l>\xi & \quad k_{3}<k_{1}<k_{4}<k_{2}\quad K\in\lbrack
k_{1},k_{4}] \notag
\end{align}
The Hamiltonian on the third reduced phase space is
\begin{eqnarray*}
&&\overline{\mathcal{H}}_{\Xi,L_{1}}=
\dfrac{3n}{4}\left(  3\beta^2-2\right)
K^{2}+\xi l(1-\beta^2)K\\
&&\hspace{1.2cm}+\dfrac{n}{2}\left(  4-\beta^2\right)  N
+n^3(\dfrac{3}{2}+\dfrac{\beta^2}{4})- \left(  {l}^{2}+{\xi}^{2}
\right) (\dfrac{\beta^2 }{2}+1)\dfrac{n}{2}
\end{eqnarray*}

Note that the reduced phase spaces as well as the Hamiltonian are  invariant under the discrete symmetry $S\rightarrow -S$. To see how to exploit it to obtain a full reducction see \cite{cushman2000}, \cite{efstathiou2004} and \cite{inarrea2004}. We choose not to further reduce our reduced phase space with respect  to these discrete symmetries, because the three dimensional picture makes it easy to access  information about the reduced orbits, and in this way one does not introduce additional critical points (fixed  points) which need special attention.

In $(K,N,S)$-space the energy surfaces are
parabolic cylinders. Notice that for $\beta^2=2/3,$ the function
$\mathcal{H}$ is linear in the variable space
$(K,N,S)$. Likewise for $\beta^2=1$, $\mathcal{H}$
modulo constants is independent
of $\xi$ and $l$. Moreover for $\beta^2=4$, $\mathcal{H}$ is only a
function of
$K$. Since 
$(V_{n\,\xi\,l},\{\cdot,\cdot\}_{3},\overline{\mathcal{H}}_{\Xi,L_{1}})$
is a Lie-Poisson system, the corresponding dynamics is given by
\begin{eqnarray}\label{sistres}
&&  \dfrac{dK}{dt}=2n(\beta^2-4)\,S,\notag\\[1.1ex]
&&  \dfrac{dN}{dt}=2[3n(3\beta^2-2)K+2\xi l(1-\beta^2)]\,S,\\[1.1ex]
&&  \dfrac{dS}{dt}=n(\beta^2-4)(K^{2}-(\xi^{2}+l^{2}+n^{2}))K\notag\\
&&\hspace{1cm}- (3\beta^2-2)[6nKN+4\xi l(\beta^2-1)N+2ln^{2}\xi].\notag
\end{eqnarray}
This system can be integrated by
means of elliptic functions, but after a classification of the
different types of flows is made, as functions of the integrals and
the parameter. Only then we will be ready for the
integration of a specific initial
value problem.  Note that in the search of relative equilibria at the
thrice reduced level, the geometry involved is really helpful; the
intersection of the Hamiltonian surfaces with the reduced phase space
gives the trajectories of the reduced system. Then, tangency of the
Hamiltonian surface with the reduced phase spaces provides relative
equilibria that generically correspond to three dimensional tori in
the original phase space.
Details on the analysis of this system, focusing on relative equilibria related to
$\beta^2=0$, $\xi=l$ and the case of physical interest $\xi=0$, are
contained in D\'{\i}az {\sl et al.} \cite{diaz2010}. The reader should take into account that $\beta^2$ is renamed in \cite{diaz2010} as $\lambda$.
\section{Symplectic charts for 4-D isotropic oscillators}
\label{sec:simplecticas}

For our study of perturbed oscillators, we switch now from Poisson to symplectic formulation. More precisely, our goal is to obtain the normal form of isotropic oscillators by Lie transforms, in the presence of symmetries.

Different charts are used for the symplectic treatment of the harmonics oscillators. Among them we find complex notation \cite{cushman1997}, action-angle variables of Poincar\'e type \cite{sanders2007}, \cite{Ferrer2000}, or variations of them in the case of resonances \cite{Deprit1999}. In our case, as we also need to reduce by the actions associated to the rotational symmetries, we have made use of a symplectic chart which incorporates the functions  (\ref{integrales}) as momenta. Moreover the connection between the isotropic and Kepler problems, two maximally superintegrable systems, allows to introduce the Delaunay transformation \cite{Deprit82}. The result is the {\sl 4-D Delaunay transformation} proposed by one of the authors \cite{Ferrer2010}, that we will use in what follows. There is an alternative procedure considering the Lissajous transformation \cite{Deprit91}, which will be presented elsewhere.

\subsection{The family of Euler projective transformations}
Let $F(\rho)$ be a smooth real function which is positive in its
domain. We consider the  family of transformations:
$\mathcal{ PE}_F: (\rho,\phi,\theta,\psi)\rightarrow (q_1,q_2,q_3,q_4)$,
dubbed as {\sl Projective Euler} variables, given by
\begin{eqnarray}\label{parametros}
\vspace{-0.3cm}
&&\hspace{-0.8cm}q_1=F({\rho})\, \sin\frac{\theta}{2}
\cos\frac{\phi-\psi}{2},\quad
q_3=F({\rho})\, \cos\frac{\theta}{2} \sin\frac{\phi+\psi}{2}, \\
&&\hspace{-0.8cm}q_2=F({\rho})\, \sin\frac{\theta}{2}
\sin\frac{\phi-\psi}{2},\quad
q_4=F({\rho})\,\cos\frac{\theta}{2} \cos\frac{\phi+\psi}{2},\nonumber
\end{eqnarray}
with $(\rho,\phi,\theta,\psi)\in R^+\times[0,2\pi)\times (0,\pi)\times
\left(-\frac{\pi}{2},\frac{\pi}{2}\right)$. For $F({\rho})=1$, this
transformation defines {\sl Euler parameters} as functions of Euler
angles.  In this paper, we only consider the case
$F(\rho)=\sqrt{\rho}$.

The canonical extension associated to the transformation
(\ref{parametros})
is readily obtained as a Mathieu transformation, which satisfies
$\sum Q_i\,dq_i= P \,d\rho + \Theta\, d\theta + \Psi\, d\psi + \Phi
\,d\phi.$
The relations among the momenta are given by
\begin{eqnarray}\label{momentos2}
&& P = \frac{1}{2\sum q_i^2}(q_1Q_1+q_2Q_2+q_3Q_3+q_4Q_4),\nonumber\\
&&\Theta=\frac{(q_1Q_1+q_2Q_2)(q_3^2 +q_4^2) -(q_3Q_3+q_4 Q_4)(q_1^2
+q_2^2)}
{2\sqrt{(q_1^2 + q_2^2)(q_3^2+q_4^2)}},\nonumber\\
&&\Psi=\frac{1}{2}(-q_2Q_1+q_1Q_2+q_4Q_3-q_3Q_4),\\
&&\Phi=\frac{1}{2}(q_2Q_1-q_1Q_2+q_4Q_3-q_3Q_4),\nonumber
\end{eqnarray}
Later on we will need the inverse transformation given by
\begin{eqnarray}\label{Eulerinversa}
&&\rho= q_1^2 +q_2^2 +q_3^2 +q_4^2,\nonumber\\ [1ex]
&&\sin\theta =\frac{2\sqrt{(q_1^2 +q_2^2)(q_3^2 +q_4^2)}}
{q_1^2 +q_2^2 +q_3^2
+q_4^2}, \qquad \cos\theta =\frac{q_3^2 +q_4^2-
q_1^2 -q_2^2 }{q_1^2 +q_2^2 +q_3^2
+q_4^2},\nonumber\\
&&\sin\psi = \frac{q_1q_3+q_2q_4}
{\sqrt{(q_1^2 +q_2^2)(q_3^2 +q_4^2)}},
\qquad \cos\psi = \frac{q_1q_4-q_2q_3}
{\sqrt{(q_1^2 +q_2^2)(q_3^2 +q_4^2)}},
\\
&&\sin\phi = \frac{q_1q_3-q_2q_4}
{\sqrt{(q_1^2 +q_2^2)(q_3^2 +q_4^2)}},
\quad\quad \cos\phi = \frac{q_1q_4+q_2q_3}
{\sqrt{(q_1^2 +q_2^2)(q_3^2 +q_4^2)}}.\nonumber
\end{eqnarray}
\begin{proposition}\label{EulerProjHam}
The isotropic oscillator $\mathcal{ H}_2$ is expressed by mean of the Euler projective variables as 
\begin{equation}
\mathcal{ H}_\omega(\rho,\theta,-,-,P,\Psi,\Theta,\Phi) =  \frac{ \omega\,\rho}{2} +2\rho P^2 +
\frac{2}{\rho} \left(\Theta^2  +\frac{\Psi^2+\Phi^2-2\,\Phi\,\Psi\,\cos\theta}{\sin^2\theta} \right).
\end{equation}
\end{proposition}
\begin{Proof} 
Considering any given initial condition for the Hamiltonian $\mathcal{ H}_\omega=h$, by applying relations (\ref{parametros}) and (\ref{momentos2}), we obtain, by straight forward computations, the above expression for the isotropic oscillator $\mathcal{ H}_2$.
\end{Proof} 

Note that the Hamiltonian has two cyclic variables, which manifests the two symmetries associated with our system we have refer in (\ref{integrales})
\subsection{From Euler to Projective Andoyer variables}
Andoyer \cite{Andoyer,Deprit1967} symplectic variables  are well known in
rotational dynamics and recently in attitude and control. Most often they
are denoted by $(\beta^2, \mu,\nu, \Lambda,M,N)$ or $(h, g, \ell, H, G, L)$.
In the following, in order to avoid confusion with invariants
notation of previous section, we propose to use
$(u_1,u_2,u_3,$ $U_1,U_2,U_3)$ for referring to them. Moreover, as in 3-D, it is convenient to
introduce the following functions
\begin{equation}\label{funcionestado}
c_1=\cos \sigma_1=  U_1/U_2, \qquad c_2=\cos \sigma_2 =  U_3/U_2,
\end{equation}
with $s_i=\sqrt{1-c_i^2}$. Hence the transformation from Euler to
Andoyer
$(\phi,\theta,\psi,$ $\Phi,\Theta,\Psi)\rightarrow(u_1,u_2,u_3,U_1,U_2,U_3)$
is given by
\begin{eqnarray}\label{proyectivaAndoyer}
&&\hspace{-1cm}\cos\theta = c_1c_2  + s_1\,s_2\,\cos u_2 ,
\hspace{1.8cm}\Phi= U_3, \nonumber\\[1.1ex]
&& \hspace{-1cm}\sin (\psi-u_1)=  \frac{\sin u_2}{\sin\theta}s_2,
\hspace{2.4cm}\Psi= U_1,\\
&&\hspace{-1cm} \sin (\phi-u_3)=  \frac{\sin u_2}{\sin\theta}\,s_1,
\hspace{2.4cm}\Theta = U_2\,\sqrt{1  - \frac{c_1^2 + c_2^2 -
2\,c_1\,c_2\,\cos\theta}{\sin^2\theta}}.\nonumber
\end{eqnarray}
with $|U_1|<U_2$ and $|U_3|<U_2$. By adding the variables $(\rho, P)$ we get a 4-D set
$(\rho,u_1,u_2,u_3, P, U_1,U_2,U_3)$ of symplectic variables that we
call `Projective Andoyer variables'.
\begin{theorem} {\sl In Projective Andoyer variables, the system
defined by (\ref{h2}), regularized by $ds=1/(4\rho)\,d\tau$, includes
the Keplerian system for any value of the integral $U_3$.}
\end{theorem}
\begin{Proof}
Following Poincar\'e technique, we first regularize  (\ref{h2}) and
our new Hamiltonian will be  $\mathcal{ K} =(\mathcal{
H}_{\omega}-h)/(4\rho)$. Considering now Proposition \ref{EulerProjHam}, after some straightforward calculations,
the Hamiltonian of the 4-D isotropic oscillator is given by
\begin{equation}\label{ESAT}
\tilde{\mathcal{ K} }= \frac{1}{2}\Big(P^2 + \frac{U_2^2}{\rho^2} \Big) -
\frac{\gamma}{\rho}
\end{equation}
in the manifold $ \tilde{\mathcal{ K}} = -  \omega/8$ and $\gamma=\frac{h}{4}$. But the above Hamiltonian corresponds to the Kepler system in polar nodal
variables, see \cite{Deprit81}. \end{Proof}

\subsection{Delaunay symplectic chart}

Following \cite{Deprit82} we plan to normalize perturbed isotropic systems by Lie-transforms  {\sl a la} Delaunay  in the following Section. Thus, we
still need to implement (in part of the phase space) another
canonical transformation of Hamilton-Jacobi type: $\mathcal{ D}_\gamma:
(L,G,\ell, g)\rightarrow (\rho, u_2,P, U_2)$. The process hinges on
the H-J equation built on the function (\ref{ESAT}). For details of
this transformation $\mathcal{ D}_\gamma$ we refer to \cite{Deprit82};
here we only include the final expressions. Among the relations
defining Delaunay transformation, which can only be given in implicit
form, we take
\begin{eqnarray}\label{varDelaunay}
&&u_2 = g + f,\hspace{2.5cm} U_2=G,\nonumber\\[1.2ex]
&& \rho= a\,(1-e\,\cos E),\hspace{1.1cm} P= \frac{L\,e\,\sin
E}{a\,(1-e\,\sin E)},\end{eqnarray}
where $a=a(L)$, $e=e(L,G)$ are given by
\begin{equation}\label{funcionesestado}
a= L^2/\gamma, \qquad \eta= G/L,\qquad e =\sqrt{1-\eta^2},\qquad
\end{equation}
and $f$ and $E$ are auxiliary angles.
The symplectic variable $\ell$ is related to them by
\begin{equation}\label{}
 \ell = E- e\,\sin E,\qquad \tan f/2 = \sqrt{(1+e)/(1-e)}\, \tan
E/2.
\end{equation}

Later on we also need\, $1/\rho= (1-e \,\cos f)/p$, where
$p=a\eta^2$, and
\begin{equation}\label{verdaderaexcentrica}
 \rho \,\cos f = a(\cos E -e), \quad   \rho \,\sin f = a\eta\,\sin E,
\quad d\ell = (1-e\,\cos E)\,dE.
\end{equation}
Completing the functions  of the momenta given above
(\ref{funcionesestado}),  it is convenient to introduce two more
state functions
\begin{equation}\label{funcionestado2}
w=U_1/L,\qquad  z=U_3/L.
\end{equation}
The previous results may be summarized in the following theorem:
\begin{theorem}
By composing the  three symplectic transformations
\medskip
\par\noindent
\hfill
\hbox{$
\begin{array}[t]{ccc}
\Big(\!\!\begin{array}{cccc}
q_1,&\!\!q_2,&\!\!q_3,&\!\!q_4\\
Q_1,&\!\!Q_2,&\!\!Q_3,&\!\!Q_4
\end{array}\!\!\Big) &
 \stackrel{\rm Projective\,\,Euler }{\longrightarrow} &
 \Big(\begin{array}{cccc}
\rho,&\!\!\phi,&\!\!\theta,&\!\!\psi\\
P,&\!\!\Phi,&\!\!\Theta,&\!\!\Psi\end{array}\Big)\\
\Biggm\downarrow{} &  & \Biggm\downarrow{{\rm P\,\!A}} \\
\Big(\begin{array}{cccc}
\ell,&\!\!u_1,&\!\!g,&\!\!u_3\\
L,&\!\!U_1,&\!\!G,&\!\!U_3
\end{array}\Big)  &
\stackrel{\rm Delaunay}{\longleftarrow}   &
\Big(\begin{array}{cccc}
\rho,&\!\!u_1,&\!\!u_2,&\!\!u_3\\
P,&\!\!U_1 ,&\!\!U_2,&\!\! U_3\end{array}\Big)
\end{array}$}
\hfill
\medskip
\par\noindent
together with the regularization ${\rm d}s=1/(4\rho)\,{\rm d}\tau$,
the  Hamiltonian of the 4-D isotropic oscillator  is given by
\begin{equation}\label{ham4Delaunay}
\mathcal{ H}_0 = -\frac{\gamma^2}{2L^2}
\end{equation}
in the manifold $\mathcal{ H}_0=-\omega/8$.
\end{theorem}
The set of variables $(\ell,g,u_1,u_3,L,G,U_1,U_3)$ is what we call
4-D Delaunay chart. As in 3-D with the classical Delaunay variables,
they represent a generalized set of action-angle variables (see
\cite{Neko72}). In the next section we show the interest of this
symplectic transformation. For the benefit of the reader let us
mention that Moser and Zehnder \cite{Moser05} introduced action-angle variables
for the Kepler problem in $\R^n$ and called them Delaunay variables.
When we restrict to $n=4$, those variables do not coincide with the
set built in this paper. \par
Finally, in order to compare results of our analysis of relative
equilibria, we need to establish the connection between the
variables $K$ (Poisson approach) and $G$ (symplectic approach)
defining the thrice reduced space. After some algebraic
manipulations, we obtain
\begin{equation}\label{conexion}
4\,G^2= \frac{1}{2}(n^2 +\xi^2+l^2) - \frac{1}{2}K^2 - N,
\end{equation}
which expresses $G$ as function of $K$ and $N$, invariants defining
the thrice reduced space, and the first integrals.

\section{Delaunay normalization of perturbed 4-D isotropic oscillators}
\label{sec:first}

Our goal in this section is to make normalization by Lie transforms \cite{Deprit69}, using 4-D Delaunay variables. As we will see  they allow to express perturbed isotropic oscillators in a `convenient' form, to implement Lie transforms up to higher order in an efficient way. What we mean by convenient is the following.

It is well known that the isotropic oscillator $\mathcal{
H}_0$, as the Kepler system, defines a maximally superintegrable
Hamiltonian system (see Fass\`o \cite{Fasso05}), whose flow is made
of periodic orbits. In that case Cushman \cite{Cushman1984} proved
that any smooth function $F$ over its phase space may be decomposed
in a unique way into a sum $F=F^{o} + F^*$ with the following
properties: {\sl (i)} $\{F^{o},\,\mathcal{ H}_0\}=0$. In other words, $F^{o}$ belongs to the kernel of the Lie derivative generated by $\mathcal{ H}_0$: $\mathcal{
L}_0: F\rightarrow \{F,\,\mathcal{ H}_0\}$; \, {\sl (ii)} There exists a
smooth function $\tilde F$ such that $\{\tilde F,\,\mathcal{ H}_0\}=F^*$;
in other words $F^*$ belongs to the image of the operator $\mathcal{ L}_0$.
In particular when the Lie transform  (see \cite{Deprit69}, \cite{Meyer09})
\[\mathcal{ L}:(\ell,g,u_1,u_3,L,G,U_1,U_3)\rightarrow
(\ell',g',u'_1,u'_3,L',G',U'_1,U'_3)\]
is carried out {\it a la} Delaunay
\cite{Deprit82}, the {\sl homological equation}, which relates the new Hamiltonian $\sum (\epsilon^j/j!)\mathcal{ H}'_j$ with the generating function $\sum(\epsilon^j/j!)\mathcal{ W}'_j$, is given by
\begin{equation}\label{homological}
\{\mathcal{ W}_j,\,\mathcal{ H}_0\}+\mathcal{ H}'_j= \tilde {\mathcal{ H}}_j, \qquad j\geq 1.
\end{equation}
Then, considering the splitting coming from Cushman theorem, the equation (\ref{homological}) may be solved choosing
\begin{equation}\label{homological2}
\mathcal{ H}'_j \,=\,\, <\! \tilde{ \mathcal{ H}}_j \!\!> = \frac{1}{2\pi} \int_0^{2\pi} \tilde{ \mathcal{ H}}_j\,{\rm d}\ell,   \qquad \mathcal{ W}_j= \frac{1}{n}\int (\tilde{ \mathcal{ H}}_j- <\! \tilde{ \mathcal{ H}}_j \!\!>)\,{\rm d}\ell,
\end{equation}
where, according with (\ref{ham4Delaunay}), $\partial\mathcal{ H}_0/\partial L=\gamma^2/L^3$.
The Hamiltonian (\ref{hamiltoniano0}) of perturbed isotropic
oscillators in Delaunay variables is given by
\[ \mathcal{ H}(\ell,g,u_1,u_3,L,G,U_1,U_3) = -\frac{\gamma^2}{2L^2} +
\varepsilon\, \mathcal{P} (\ell,g,u_1,u_3,L,G,U_1,U_3;\beta).\]
and the normalized Hamiltonian up to order $k$ takes the form
\[ \mathcal{ H}' = -\frac{\gamma^2}{2L'^2}+ \sum^k
\frac{\varepsilon^n}{j!} \mathcal{ H}'_j
(-,g',u'_1,u'_3,L',G',U'_1,U'_3;\beta) + \mathcal{
O}(\varepsilon^{k+1}).\]
The main feature of the process is that, at each order, solving the
{\sl homological equation} (\ref{homological}) only involves quadratures. \par

\subsection{Delaunay normalization and symmetries}
Apart from the procedure associated to Delaunay normalization
we have just mentioned, when we restrict to systems with
the symmetries (\ref{integrales}), the use of Delaunay chart
shows its full advantage. Indeed, in that case  $u_1$ and $u_3$ are
cyclic, in other words, we have
\begin{equation}\label{hamDelaunay}
\mathcal{ H}(\ell,g,-,-,L,G,U_1,U_3) = -\frac{\gamma^2}{2L^2} +
\varepsilon\, \mathcal{P} (\ell,g,-,-,L,G,U_1,U_3;\beta;\varepsilon).
\end{equation}
In geometric language, the use of Delaunay variables has carried out
the toral reduction associated to the actions defined by the
symmetries (\ref{integrales}), and the function (\ref{hamDelaunay})
is the reduced Hamiltonian. The Hamiltonian system given by
(\ref{hamDelaunay}) is a 2-DOF system. Thus, the normalized system will take the form
\begin{equation}\label{hamDelaunay2}
\mathcal{ H}' = -\frac{\gamma^2}{2L'^2} + \sum^k
\frac{\varepsilon^n}{j!} \mathcal{ H}'_j
(-,g',-,-,L',G',U_1',U_3';\beta) + \mathcal{ O}(\varepsilon^{k+1}).
\end{equation}

In our case normalizing up to the order needed, after truncating, we will
obtain an integrable 1-DOF Hamiltonian system. From now on we drop
primes in the variables in order to simplify the expressions.

\subsection{Carrying out the normalization. The van der Waals model
as a benchmark}\label{CarryingNormalization}

In order  to implement these normalizations, since we are dealing with
polynomial perturbations, from experience with 3-DOF Kepler
systems we introduce the auxiliary variable $E$ by using the following proposition.
\begin{proposition}\label{propo1}
The functions  $F(\rho,\theta)= \rho^m\,\cos^n\theta$, $(m,n\in
\mathbb{N}$ and $m\geq n$), expressed  in Delaunay variables, could be written as
$F(\rho,\theta)=\sum_{i=0}^{m} (C_i\cos i\,E + S_i\sin i\,E),$ where $C_i$
and $S_i$ are given by expressions of the form 
$\sum_j (c_{ij}\,\cos j\,g + s_{ij}\,\sin j\,g )$ which 
belong to the kernel of the Keplerian Lie derivative,
and $c_{ij}$, $s_{ij}$ are rational functions of the momenta.
\end{proposition}
\begin{Proof}
It is straightforward, based on relations given in (\ref{proyectivaAndoyer}), (\ref{varDelaunay}) and (\ref{verdaderaexcentrica}). 
\end{Proof}

Applying the above result, the Hamiltonian (\ref{Ham1}) takes the form
\begin{equation}
\mathcal{ H} = -\frac{\gamma^2}{2L^2} +  \sum^k_{j=0}
\frac{\varepsilon^j}{j!}  \sum_{i}(C_{ij}\cos i E +S_{ij} \sin i E)+ \mathcal{ O}(\varepsilon^{k+1}),
\end{equation}
where functions $C_{ij}= C_{ij}(g, L,G,U_1,U_3;\beta)$ and $S_{ij}= S_{ij}(g,
L,G,U_1,U_3;\beta)$ belong to the kernel of the Lie transform.\par

Thus, by considering the differential relation (\ref{verdaderaexcentrica}) and by replacing in (\ref{homological2}), we obtain the part in the kernel as follows
\begin{eqnarray*}
&&<\!\mathcal{ H}_1\!\!>\,= \frac{1}{2\pi} \int_0^{2\pi} \mathcal{ H}_1\,d\ell\\
&&\hspace{1.1cm}
= \frac{1}{2\pi} \int_0^{2\pi}\sum_{i=0} (C_{i1}\cos i E +S_{i1} \sin i E)
(1-e\,\cos E)\,dE\\[1.5ex]
&&\hspace{1.1cm}=C_{01} + C_{11} \cos g + C_{21} \cos 2g.
\end{eqnarray*}
with $C_{i1} = C_{i1} (L,G,U_1,U_3;\beta)$ given by
\begin{eqnarray}\label{coef}
&& C_{01} = \frac{1}{16} a^2 (2 + 3e^2) [\,2 +  (\beta^2-1)(2 c_1^2
c_2^2 + s_1^2 s_2^2)\,],\nonumber\\[1ex]
&&C_{11} = -\frac{1}{4} a^2  (\beta^2-1)(4 + e^2)\, e\,c_1 c_2 s_1 s_2,\\[1ex]
&&C_{21} = \frac{5}{16} a^2  (\beta^2-1) e^2 s_1^2 s_2^2.\nonumber
\end{eqnarray}
Then, we solve the first order homological equation (\ref{homological}). The first order normalized Hamiltonian takes the form

\begin{eqnarray*}
&&\mathcal{ H} = -\frac{\gamma^2}{2L^2}+ \varepsilon\,(C_{01} + C_{11} \cos g + C_{21} \cos 2g) + \mathcal{ O}(\varepsilon^2)
\end{eqnarray*}

The first order  $\mathcal{ W}_1$ of the generating function is obtained computing the second quadrature in (\ref{homological2}), then, taking into account the shorthand $\alpha=\beta^2-1$, we get the following expression
\begin{eqnarray*}
\label{ w1}
&&\mathcal{ W}_1=\dfrac{a^3 L}{96 \gamma } \left(\alpha(2 c_1^2  c_2 ^2+s^2 s_2 ^2)+2\right)\\
&&   \hspace{1cm}
\Big(c_1\big(e^3\sin (3 u) -9 e^2\sin (2 E) +3e(3 \eta^2+5) \sin (u)\big)+ \alpha s_1s_2 \\
&&   \hspace{1cm}
+s_1 s_2(1-\eta)^2 \big( -15 e \sin (2 g-E)   +(9+6 \eta) \sin (2 g-2 E)  \\
&&   \hspace{3.8cm}
-e \sin (2 g-3 E)\big)\\
&&   \hspace{1cm}
-s_1 s_2(1+\eta)^2 \big(15 e \sin (2 g+E)-(9-6 \eta ) \sin (2 g+2 E)\\
&&   \hspace{3.8cm}
+e \sin (2 g+3 E)\big)\\
&&   \hspace{1cm}
+4 c_1 c_2(1-\eta ) \big(e^2\sin (g-3 E) -3 e(3+\eta) \sin (g-2 E)\\
&&   \hspace{3.8cm}
+3(2 \eta ^3+\eta ^2+5) \sin (g-E)\big)\\
&&   \hspace{1cm}
-4 c_1c_2 (1+\eta) \big(e^2\sin (g+3 E) -3 e(3-\eta ) \sin (g+2 E) \\
&&   \hspace{3.8cm}
+3(-2 \eta ^3+\eta ^2+5) \sin (g+E)\big)\Big).
\end{eqnarray*}

Analogously we obtain the second order normalization for (\ref{hamDelaunay2}),  which will be used in the study of stability. The second order coefficients are

\begin{eqnarray*}
&&C_{02}=-60\,{\alpha}^{2}{c_2}^{2}c_1^{2} s_1^2s_2^2\,{e}^{6}\\
&&   \hspace{1cm}+
\Big(\frac{{\alpha}^{2}}{8}\big( 1139\,{c_2}^{4}c_1^{4}-1486({c_2}^{4}c_1^{2}+{c_2}^{2}c_1^{4})+579(c_1^{4}+{c_2}^{4})\\
&&   \hspace{2.1cm}+1268\,{c_2}^{2}c_1^{2} -518(c_1^{2}+2\,{c_2}^{2})-61\big) +63\left(\alpha \, \Lambda_{3,1} -1\right) \Big){e}^{4}\\
&&   \hspace{1cm}
+\Big( {\alpha}^{2} \,\big( 2089\,{c_2}^{4}c_1^{4}-1756({c_2}^{2}c_1^{4}+{c_2}^{4}c_1^{2})+399(c_1^{4}+{c_2}^{4}) \\
&&   \hspace{2.1cm}+2048\,{c_2}^{2}c_1^{2}-628(c_1^{2}+{c_2}^{2})+229\big)-396\left( \alpha\, \Lambda_{3,1} -1\right)\Big){e}^{2}\\
&&   \hspace{1cm}
+\Big( {\alpha}^{2}\,\big( 557\,{c_2}^{4}c_1^{4}-382({c_2}^{4}c_1^{2}+382{c_2}^{2}c_1^{4})+17(c_1^{4}+{c_2}^{4})\\
&&   \hspace{2.1cm}+308\,{c_2}^{2}c_1^{2}-22({c_2}^{2}+c_1^{2})+5\big)- 96\left(\alpha\, \Lambda_{3,1} -1\right)\Big),
 \end{eqnarray*}
\begin{eqnarray*}
 &&C_{12}=c_2\,c_1s_1s_2 \,\Big( \big(-155{\alpha}^{2}\small( \Lambda_{2,1}+41c_1^2c_2^2+36\small) +192\,\big){e}^{5}\\
&&   \hspace{2.85cm}
-(447{\alpha}^{2} \left( \Lambda_{2,1}-19c_1^2c_2^2-290 \right) +708\alpha\,+2){e}^{3}\\
&&   \hspace{2.85cm}
-\,\big(67\, {\alpha}^{2}\,\left( \Lambda_{2,1}+c_1^2c_2^2+27 \right)-1200\,\alpha-16\,\big)e\Big),\\\\
 &&C_{22}=s_1s_2\Big(-60\,{\alpha}^{2}{c_2}^{2}c_1^{2}  {e}^{6}\\
&&   \hspace{2cm}
+ \big(\frac{{-\alpha}^{2}}{2} \left( 71\Lambda_{5,2}+8c_1^2c_2^2-14 \right)-57\,\alpha\,
 \big) {e}^{4}\\
&&   \hspace{2cm}
+ \big({\alpha}^{2}\left( -147\Lambda_{2,1}+25c_1^2c_2^2-66 \right) +486\,\alpha\,+3\big)e^2  \Big),\\\\
 &&C_{32}=\,c_1s_{1}^3c_2 s_{2}^3{\alpha}^{2}\Big(85{e}^{5}+10{e}^{3}\Big),\\\\
  &&C_{42}=-{\frac {95}{8}}{\alpha}^{2}s_1^2s_2^2\,{e}^{4},
 \end{eqnarray*}
 where
 \begin{equation}
\label{ }
\Lambda_{n,i}=c_1^2+c_2^2(-1)^i n\,c_1^2c_2^2(-1)^i.\nonumber
\end{equation}

As the van der Waals system includes several integrable cases, in order to study stability in those cases, higher order normalization, at least second order, is needed. This subject will be tackled in future works.

Here we try only to identify {\sl relative equilibria}. In fact, this
is the procedure to follow in order to classify
and compute different types of orbits defined by our system. For the generic case, first order normalization is sufficient. The second order normalization is required when we consider degenerate situations, connected with integrable cases. For this reason we give both, first and second order in Section \ref{CarryingNormalization}.

\section{Searching for relative equilibria and invariant tori of the
first order normalized system}
\label{sec:periodic}

Focusing on the task that we have announced before, we use the first order normalization of the Hamiltonian. Thus the corresponding differential system is given by
\begin{eqnarray}
&&\dot \ell = \frac{\partial \mathcal{ H}}{\partial L} =
\frac{\gamma^2}{L^3} + \varepsilon \frac{\partial \mathcal{P}}{\partial
L} , \quad\qquad \dot L = -\frac{\partial \mathcal{ H}}{\partial \ell}
=0,\label{sn1}\\
&&\dot g = \frac{\partial \mathcal{ H}}{\partial G} = \phantom{ \omega L
+}\,\varepsilon \frac{\partial \mathcal{P}}{\partial G} ,\quad\qquad
\dot G =- \frac{\partial \mathcal{ H}}{\partial g} ,\label{sn2}\\
&&\dot u_1 = \frac{\partial \mathcal{ H}}{\partial U_1} = \phantom{
\omega L +}\varepsilon
\frac{\partial \mathcal{P}}{\partial U_1},\qquad \dot U_1=-
\frac{\partial \mathcal{ H}}{\partial u_1} =0 \label{sn3}\\
&&\dot u_3 = \frac{\partial \mathcal{ H}}{\partial U_3} = \phantom{
\omega L +}\varepsilon \frac{\partial \mathcal{P}}
{\partial U_3} ,\qquad \dot U_3=- \frac{\partial \mathcal{ H}}{\partial
u_3} =0\label{sn4}
\end{eqnarray}
In other words $L,U_1$ and $U_3$ are integrals, as we already know,
and the system splits into a 1-DOF Hamiltonian system,
namely
\begin{equation}\label{3reducido}
\dot g = \varepsilon \frac{\partial \mathcal{P}}{\partial G},\qquad
\dot G =-\varepsilon \frac{\partial \mathcal{P}}{\partial g}
\end{equation}
and three quadratures coming from
\begin{equation*}
\dot \ell = \frac{\gamma^2}{L^3}+ \varepsilon \frac{\partial \mathcal{
P}}{\partial L} , \qquad \dot u_1 = \varepsilon
\frac{\partial \mathcal{P}}{\partial U_1}, \qquad \dot u_3 = \varepsilon
\frac{\partial \mathcal{P}}{\partial U_3},
\end{equation*}
as soon as the system Eqs.~(\ref{3reducido}) is solved. We are not
going to explore this way, which moreover
involves hyperelliptic integrals. In other words, it will be more
convenient to rely on numerical integrations.

As we have already said, our previous work on this model  concentrated on searching for
relative equilibria related to singular points of the energy-moment map.
We have already mentioned that the complete analysis should be done
using invariants. In fact the lower dimensional relative equilibria, i.e. those that correspond to invariant $\Bbb{S}^1$ or $\T^2$, are given by the
singularity of the moment map for the $\T^3$-symmetry group, and can
be described by a moment polytope; for details see \cite{diaz2010}.\par

Excluding those solutions, in the remaining open domain we may use symplectic charts. Relative equilibria of the system (\ref{sn1}) - (\ref{sn4}) may be classified as follows:
\par\noindent
$\bullet$ {\bf invariant 3-tori}. They are
the solutions $(g,G)$ of the  system defined by
\begin{equation}\label{sistemareducido3}
\dot G =- \varepsilon\frac{\partial \mathcal{P}}{\partial g}=0 ,\qquad
\dot g = \varepsilon \frac{\partial \mathcal{P}}{\partial G}=0.
\end{equation}
as functions of $U_1$, $U_3$ and $\beta$.\par\noindent
$\bullet$ {\bf invariant 2-tori}. In a similar way, the search for
{\sl invariant 2-tori} splits into the search the roots of the subsystems. Namely
\begin{equation}\label{sistemareducido2}
\dot G =- \varepsilon\frac{\partial \mathcal{P}}{\partial g}=0 ,\qquad
\dot g = \varepsilon \frac{\partial \mathcal{P}}{\partial G}=0,  \qquad
\dot u_1 = \varepsilon \frac{\partial \mathcal{P}}{\partial U_1}=0,
\end{equation}
where the possible roots will be function of $U_3$ and $\beta$,
and
\begin{equation}\label{sistemareducido22}
\dot G =- \varepsilon\frac{\partial \mathcal{P}}{\partial g}=0 ,\qquad
\dot g = \varepsilon \frac{\partial \mathcal{P}}{\partial G}=0,  \qquad
\dot u_3 = \varepsilon \frac{\partial \mathcal{P}}{\partial U_3}=0,
\end{equation}
where the possible roots will be function of $U_1$ and $\beta$.
\par\noindent
$\bullet$ {\bf periodic orbits}: The research for {\sl periodic
orbits} of the original system is equivalent \cite{Yanguas2008}
to find the roots of the first reduced system of equations
\begin{equation}\label{sistemareducido}
\dot G =- \varepsilon\frac{\partial \mathcal{P}}{\partial g}=0,\quad
\dot g = \varepsilon \frac{\partial \mathcal{P}}{\partial G}=0, \quad
\dot u_1 = \varepsilon
\frac{\partial \mathcal{P}}{\partial U_1}=0, \quad \dot u_3 =
\varepsilon \frac{\partial \mathcal{P}}{\partial U_3}=0 ,
\end{equation}
Solutions of the system given by Eqs.~(\ref{sistemareducido3})   define  invariant 3-tori $\T^3(\ell,u_1,u_3)$, except for the singular points \cite{diaz2010}. The system defined by Eqs.~(\ref{sistemareducido2}) gives the invariant 2-tori $\T^2(\ell,u_3)$, each of them is attached to a 2-torus $\T^2(g,u_1)$ made of fixed points. Likewise  the system defined by Eqs.  (\ref{sistemareducido22}) is the invariant 2-tori $\T^2(\ell,u_1)$, each of them is attached to a 2-torus $\T^2(g,u_3)$ made of fixed points. Finally, the intersection of the previous types of 2-tori, defined by the system of Eqs.~(\ref{sistemareducido}),  correspond to periodic solutions
$S^1(\ell)$, each of them is attached to a 3-torus $\T^3(g,u_1,u_3)$ made of fixed points. They correspond to the short period solutions related to the unperturbed system. In what follows we search for invariant 3-tori and periodic solutions of short period.

\subsection{Searching for invariant 3-tori}

To find families of relative equilibria (periodic orbits, invariant tori, etc), we start searching for invariant 3-tori, the common condition for all the cases. Explicitly, according to above paragraphs, the equations (\ref{sistemareducido3}) are
\begin{eqnarray}
&&\frac{\partial C_{01} }{\partial G} + \frac{\partial C_{11} }{\partial G}\cos g + \frac{\partial C_{21} }{\partial G}\cos 2g =0 ,\label{sistemareducido31}\\ [1.5ex]
&&(C_{11}  + 4\,C_{21} \cos g)\,\sin g=0.\label{sistemareducido32}
\end{eqnarray}
Due to the structure of Eq. (\ref{sistemareducido32}), and keeping in mind the domain of existence of 4-D Delaunay variables, the only possibility leading to roots is  $\sin g=0$.  Replacing  in  (\ref{sistemareducido31}),  we have an equation to be solved for each of the values of $\cos g$, namely $\pm 1$. Such equations are:
\begin{alignat}{2}\label{equacionesSeng0}
\frac{\partial C_{01} }{\partial G} + \frac{\partial C_{11} }{\partial G} + \frac{\partial C_{21} }{\partial G} =0, & \qquad
\frac{\partial C_{01} }{\partial G} - \frac{\partial C_{11} }{\partial G} + \frac{\partial C_{21} }{\partial G} =0.
\end{alignat}
Both equations can be solved at one time applying the following strategy. By using the state functions   (\ref{funcionesestado})  and (\ref{funcionestado2}), both equations (\ref{equacionesSeng0}) can be rewritten
in the form;
$R(\eta,z,w)x(\eta,z,w)+Q(\eta,z,w)=0\quad \text{and} \quad R(\eta,z,w)x(\eta,z,w)-Q(\eta,z,w)=0$, in which:
\begin{eqnarray}
&&x(\eta,z,w)=\sqrt{(1-w^2)(1-z^2)e} ,\nonumber\\ [1ex]
&&R(\eta,z,w)=(-(3+4\alpha)\eta^6+(5w^2+7w^2z^2+5z^2)\,\alpha\,\eta^2-20w^2z^2\alpha)\,\eta^2 ,\nonumber\\ [1ex]
&&Q(\eta,z,w)=\alpha \,w\,z\left[\eta^8 -(10+11w^2+w^2z^2)\eta^4 -11z^3 \eta^3\right. \\
&&\hspace{2cm}\left.+(15w^2+17w^2z^2 +15z^2)\eta^2 +5z^5\,\eta -20w^2z^2\right]\, .\nonumber
\end{eqnarray}
By multipliying equations $Rx+Q=0$ and $Rx-Q=0$ we obtain a new polynomial
$P(\eta)=R^2x^2-Q^2$ in $\eta$ as independient variable, where $w$, $z$ and $\alpha$ are taken as parameters. The roots of $P$ represent the set of solutions
for the original equations (\ref{equacionesSeng0}). Explicitly  the polynomial is:
\begin{equation}\label{polinomio}
P(\eta)=a_{6}{\eta}^{6}+a_{5}{\eta}^{5}+a_{4}{\eta}^{4}+a_{3}{\eta}^{3}+a_{2}{\eta}^{2}+a_{1}{\eta}+a_{0},
\end{equation}
where the coefficients are given by
\begin{align*}
a_{6}=& -9-24\,\alpha-16\,{\alpha}^{2},\\
a_{5}=& \,{w}^{2}(16\,{\alpha}^{2}+9-{\alpha}^{2}{z}^{2}+24\alpha)+9+24\,
\alpha+16\,{\alpha}^{2}+{z}^{2}(24\,\alpha+9\,+16\,{\alpha}^{2}),\\
a_{4}=&\,{w}^{2} ( 24\,{\alpha}^{2}+18\,\alpha\,{z}^{2}+6\,\alpha\,-9\,-9\,{z}^{2}
+30\,{\alpha}^{2}{z}^{2}{w}^{2})+{z}^{2}(6\,\alpha\,+24\,{\alpha}^{2}-9),\\
a_{3}=&\, {w}^{4}( -40\,{\alpha}^{2}-42\,
\alpha\,{z}^{2}-30\,\alpha\,-34\,{\alpha}^{2}{z}^{2}+2\,{z}^{4})+9\,{z}^{2}{w}^{2}-30\,{z}^{2}-40\,{\alpha}^{2}{z}^{2}\\
&-40\,{\alpha}^{2}{z}^{4}+{\alpha}{w}^{2}(198\,{z}^{2}-34{\alpha}{z}^{4}+30\,+285\,{\alpha}{z}^{2}+40\,{\alpha}+
42\,{z}^{4})-30\,\alpha\,{z}^{4}, \\
a_{2}=& \,\alpha\, ( 30\,{w}^{4}+42\,{z}^{4}{w}^{4}+192\,
{z}^{2}{w}^{4}+15\,\alpha\,{w}^{4}-17\,\alpha\,{z}^{4}{w}^{4}+266\,
\alpha\,{z}^{2}{w}^{4}-30\,{z}^{4}  \\
&+266\,\alpha\,{z}^{4}{w}^{2}+192\,{z}^{4}{w}^{2}
+180\,{z}^{2}{w}^{2}+290\,\alpha\,{z}^{2}{w}^{2}+15\,\alpha\,{z}^{4}),   \\
a_{1}=&\,\alpha^2\, (25\,{w}^{4}-{w}^{6}{z}^{6}+27\,{w
}^{6}{z}^{4}-51\,{w}^{6}{z}^{2}+25\,{w}^{6}+27\,{z}^{6}{w}^{4}-139\,
{z}^{4}{w}^{4}\\
&-225\,{z}^{2}{w}^{4}+25\,{z}^{6}+25\,
{z}^{4}-51\,{z}^{6}{w}^{2}
-225\,{z}^{4}{w}^{2}-50\,{z}^{2}{w}^{2})\\
&+150\,{z}^{2}{w}^{4}+150\,{z}^{4}{w}^{2}+162\,{z}^{4}{w}^{4},\\
a_{0}=&\,5\,\alpha\, (\alpha\,{w}^{6}{z}^{4}-3\,\alpha\,{w}^{6}{z}^{6}-
5\,\alpha\,{w}^{6}+7\,\alpha\,{w}^{6}{z}^{2}+24\,{z}^{4}{w}^{4}+\alpha
\,{z}^{6}{w}^{4} \\
&+5\,\alpha\,{z}^{2}{w}^{4}+18\,\alpha\,{z}^{4}{w}^{4}-
5\,\alpha\,{z}^{6}+5\,\alpha\,{z}^{4}{w}^{2}+7\,\alpha\,{z}^{6}{w}^{2}).
\end{align*}
This polynomial corresponds with the one obtained from the system (\ref{sistres}) when we look for relative equilibria within the Poisson approach \cite{diaz2010}. The complexity in both analyses is similar because the polynomials are of the same degree. 
It is straighforward to check that when we impose the constraint $z=0$, we recover the expressions of the study done with the classical Delaunay variables \cite{ElipeFerrer94}.

The search for possible invariant $\T^2$ solutions of Systems (\ref{sistemareducido2}) and (\ref{sistemareducido22}) requires a similar analysis. We will not refer to them here. Instead, we focus on periodic orbits, where we obtain a benefit working in symplectic formalism.

\subsection{Searching for periodic orbits}

With Poisson formalism, finding periodic solutions is in correspondence with computing relative equilibria of the system defined with the normalized Hamiltonian (\ref{primerreducido}). This study requires to deal with the constraints defining $\mathbb{CP}^{3}$ and leads to a large polynomial system; for details see \cite{diaz2010}. Although partial results are obtained, the presence of the physical parameter introduces complicated expressions.\par

The situation changes dramatically when we switch to symplectic formalism. The complete system (\ref{sistemareducido}) related to periodic solutions takes the form
\begin{equation*}\label{sistemaperiodicas}
\begin{aligned}
&&\frac{\partial C_{01} }{\partial G} + \frac{\partial C_{11} }{\partial G}\cos g + \frac{\partial C_{21} }{\partial G}\cos 2g =0, \\
&&(C_{11}  + 4\,C_{21} \cos g)\,\sin g=0,\\
&&\frac{\partial C_{01} }{\partial U_1} + \frac{\partial C_{11} }{\partial U_1}\cos g + \frac{\partial C_{21} }{\partial U_1}\cos 2g =0, \\
&&\frac{\partial C_{01} }{\partial U_3} + \frac{\partial C_{11} }{\partial U_3}\cos g + \frac{\partial C_{21} }{\partial U_3}\cos 2g =0.
\end{aligned}
\end{equation*}
The solutions of the system are the set of the periodic orbits in the first order normalized problem.  By solving the above system imposing $\alpha\neq 0$, we obtain that periodic orbits are characterized by $U_1=U_3=0$ and 
\begin{equation}\label{RaicesOrbitasPeriodicas }
\begin{aligned}
&&|U_1|=|U_3|= \sqrt {{\frac {e \cos g \left( 4+5\,e \cos g+{e}^{2} \right) +1-{e}^{2}}{e \cos g \left( 
5\, e \cos g +8+2\,{e}^{2} \right) +3+2\,{e}^{2}}}} \\
&&|U_1|=|U_3|= \sqrt {{\frac {e \cos g \left( -4+5\,e \cos g-{e}^{2} \right) +1-{e}^{2}}{e \cos g \left( 
5\,e \cos g -8-2\,{e}^{2} \right) +3+2\,{e}^{2}}}} 
\end{aligned}
\end{equation}

Then (\ref{sistemareducido}) becomes a system of three  equations  given by
\begin{equation}\label{sistemaperiodicas}
\begin{aligned}
&&\frac{\partial C_{01} }{\partial G} + \frac{\partial C_{11}
}{\partial G}\cos g + \frac{\partial C_{21} }{\partial G}\cos 2g =0, \\
&&(C_{11}  + 4\,C_{21} \cos g)\,\sin g=0,\\
&&\frac{\partial C_{01} }{\partial U_3}+ \frac{\partial C_{11}
}{\partial U_3} \,\cos g + \frac{\partial C_{21} }{\partial U_3}\,\cos 2g  =0.
\end{aligned}
\end{equation}
Again we look for solutions when $\sin g=0$. Then, we drop the second equation and we have two systems to solve related to $\cos g=\pm 1$:
\begin{equation}
(i) \qquad \frac{\partial C_{01} }{\partial G} + \frac{\partial C_{11}
}{\partial G} + \frac{\partial C_{21} }{\partial G} =0, \qquad
\frac{\partial C_{01} }{\partial U_3} + \frac{\partial C_{11}
}{\partial U_3} + \frac{\partial C_{21} }{\partial U_3} =0,
\end{equation}
and
\begin{equation}
(ii) \qquad \frac{\partial C_{01} }{\partial G} - \frac{\partial C_{11}
}{\partial G} + \frac{\partial C_{21} }{\partial G} =0, \qquad
\frac{\partial C_{01} }{\partial U_3} - \frac{\partial C_{11}
}{\partial U_3} + \frac{\partial C_{21} }{\partial U_3} =0.
\end{equation}
Considering case $(i)$, from the second equation we obtain
\[c_2\,({e}^{2}+3\,e+1 - c_2^2\,\left( 2\,e+3 \right)  \left( e+1
\right)) =0. \]
Thus, we have either $c_2=0$ or
\begin{equation}
c_2^2={\frac {1+3\,e+{e}^{2}}{ \left(3+ 2\,e \right)\left(1+e
\right) }}.
\end{equation}
If $c_2=0$, we obtain $\alpha=-3/4$. Otherwise, from the first equation we find
\begin{equation}
\alpha=\,{\frac {3e \left(3+ 2\,e \right)^{2}}{3\,{e}^{4}+2\,{e}^{3}-10\,{e}
^{2}-3\,e+8}}
\end{equation}
Likewise, for the case $(ii)$ the equation is
\[c_2\,({e}^{2}-3\,e+1 - c_2^2\,\left( 2\,e-3 \right)  \left( e-1
\right)) =0. \]
Hence, again  $c_2=0$ or
\begin{equation}
c_2^2={\frac {1-3\,e+{e}^{2}}{ \left( 2\,e-3 \right)\left(e-1
\right) }}.
\end{equation}
If $c_2=0$, we obtain $\alpha=-3/4$. Otherwise, from the first equation now we find
\begin{equation}
\alpha=-{\frac {3e \left(3 -2\,e\right) ^{2}}{3\,{e}^{4}-2\,{e}^{3}-10\,{e}
^{2}+3\,e+8}}.
\end{equation}
Thus, for each member of the Van der Waals family defined by $\alpha$, we have the relation between integrals which give periodic orbits. The period is computed replacing those values in the right hand terms of  (\ref{sn1}).
\section*{Conclusions and future work}

In contrast to some claims in \cite{cush2005}, Poisson and symplectic formulations of perturbed 4-D isotropic oscillators systems behave as complementary in the analysis of the main features of those systems. The first approach, built on the invariants and their algebraic relations, is necessary in order to carry out regular and singular reductions of the systems, in particular when singular points are present in the reduced spaces. Nevertheless the use of those invariants and their relations, introduce large computations in the search of relative equilibria, specially when physical parameters are involved. Only those equilibria related to the symmetry groups of the system are more easily computable. \par

The situation is rather different when the symplectic approach is followed. From the geometric mechanics point of view this helps to portrait the toral structure of the phase space, necessary when stability KAM theory is applied. Moreover, from the algebraic perspective, the number of equations and relations involved falls sharply because there are no constraints. Normal forms are computed very efficiently in this frame.  Of course, this is at the expense of considering only an open domain of the phase space, where those symplectic variables are defined. For this reason we should  begin with Poisson formalism. It is only after studding  singular points that we can switch to symplectic techniques with properly chosen variables.

As an illustration, the Van der Waals family is studied in detail here, showing the pros and cons of both approaches. The reconstruction process connecting the whole analysis of relative equilibria with the original system, is still to be tackled. In particular, the analysis of the different types periodic orbits.

An open question is whether there are other invariants which can lead to similar equations to those appearing in symplectic variables . In the same vein, connected with our choice as 4D-Delaunay variables [23], it might be worth considering other recent proposals [17], [30]. Results on this, now in  progress, will be published in [3].

\section*{Acknowledgements}
We deeply appreciate suggestions and comments of the referees which contribute to the improvement of the text. Support from the Government of Spain is recognized. It came in the
form of research projects MTM2008-04699-C03-03  (G.D.), MTM 2009-10767 (S.F.), ESP 2007-64068 and AYA 2009-11896 (M.L) of the Ministry of Science, and a grant 12006/PI/09 from Fundaci\'on S\'eneca of the Autonomous Region of Murcia.

\bibliographystyle{elsart-num-sort}
\bibliography{Bibliography}

\end{document}